\documentclass[a4paper]{amsart}

\usepackage[utf8]{inputenc}
\usepackage{mathtools}
\usepackage{amssymb,amsfonts,amsmath}
\usepackage{enumerate}
\usepackage{color}
\usepackage{mathrsfs}
\usepackage{graphicx}
\usepackage{verbatim}
\usepackage{tikz-cd}

\newtheorem{theorem}{Theorem}

\theoremstyle{definition}

\newtheorem*{introdefinition}{Definition}

\theoremstyle{remark}

\newcommand{\N}{\mathbb{N}}
\newcommand{\Q}{\mathbb{Q}}

\newcommand{\Z}{\mathbb{Z}}

\DeclareMathOperator{\GL}{GL}

\DeclareMathOperator{\PSL}{PSL}

\DeclareMathOperator{\SL}{SL}

\DeclareMathOperator{\T}{T}

\newcommand{\defeq}{\mathrel{\mathop{:}}=}

\renewcommand{\epsilon}{\varepsilon}

\title[Branch groups are not boundedly generated]{On representations of direct products and the bounded generation property of branch groups}
\author[S. Kionke]{Steffen Kionke}
\author[E. Schesler]{Eduard Schesler}
\address{FernUniversit\"at in Hagen \\ Fakult\"at f\"ur Mathematik und Informatik \\
58084 Hagen}
\email{steffen.kionke@fernuni-hagen.de}
\email{eduard.schesler@fernuni-hagen.de}
\thanks{Funded by the Deutsche Forschungsgemeinschaft (DFG, German Research Foundation) - 441848266}
\subjclass[2010]{Primary 20E08; Secondary 20E26}
\keywords{bounded generation, branch group, faithful representations of direct products}

\begin{document}
\begin{abstract}
We prove that the minimal representation dimension of a direct product $G$ of non-abelian groups $G_1,\ldots,G_n$ is bounded below by $n+1$ and thereby answer a question of Ab\'{e}rt.
If each $G_i$ is moreover non-solvable, then this lower bound can be improved to be $2n$.
By combining this with results of Pyber, Segal, and Shusterman on the structure of boundedly generated groups we show that branch groups cannot be boundedly generated.
\end{abstract}
\maketitle

\section*{Introduction}

An infinite group $G$ is called \emph{just-infinite} if all of its proper quotients are finite.
Obvious examples of just-infinite groups are virtually simple groups.
Other examples arise from irreducible lattices in higher rank semisimple Lie groups, such as $\SL_n(\Z)$ for $n \geq 3$, after dividing out their centers, see~\cite[Chapter IV]{Margulis91}.
Such groups are in fact \emph{hereditarily just-infinite}, which means that they are residually finite and all of their finite index subgroups are just-infinite.
Grigorchuk's group~\cite{Grigorchuk80} provided the first example of a just-infinite group that is \emph{not} virtually a finite direct power of a simple or a hereditarily just infinite group.
Grigorchuk's group is a just-infinite \emph{branch group}, which means that its commensurability classes of subnormal subgroups form a lattice that is isomorphic to the lattice of open and closed subsets of a Cantor set.
By Wilson's classification~\cite{Wilson71} just-infinite groups fall into three classes.
Every just-infinite group $G$ is either a branch group or virtually a direct power of a simple or a hereditarily just-infinite group.

%see~\ref{def:branch-group-abstract} for a precise definition.

Since its introduction by McCarthy~\cite{McCarthy68} in the late 1960's, the class of just-infinite groups remained an active field of research.
One reason might be that every finitely generated infinite group admits a just-infinite quotient.
Thus whenever there is some finitely generated, infinite group $G$ that admits a property $\mathcal{P}$ that is preserved under homomorphic images, then there is also a finitely generated just-infinite group with $\mathcal{P}$.
Following~\cite{BGS-branch}, we call a property $\mathcal{P}$ that is preserved under homomorphic images an $\mathcal{H}$-property.
Well-known examples of $\mathcal{H}$-properties include amenability, property $(\T)$, bounded generation, being a torsion group, having subexponential growth etc.
In view of Wilson's classification, it is natural to investigate which of the three classes of just-infinite groups contain groups that satisfy a given $\mathcal{H}$-property $\mathcal{P}$.
For the $\mathcal{H}$-property ``being a torsion group'' this question is settled.
In this case it is know that there are finitely generated simple groups~\cite{Olshanskii79}, just-infinite branch groups~\cite{Grigorchuk80}, and hereditarily just-infinite  groups~\cite{ErshovJaikin10} that are torsion.
On the other hand, there are torsion-free, finitely generated, just-infinite groups that are simple~\cite{BurgerMozes00}, branch~\cite{BartholdiGrigorchuk99}, and hereditarily just-infinite (e.g.\ $\Z$).

The purpose of this note is to study this question for the bounded generation property.
Recall that a group $G$ is \emph{boundedly generated} if it contains a finite subset $\{g_1,\ldots,g_n\}$ such that every $g \in G$ can be written as $g = g_1^{k_1} \cdots g_n^{k_n}$ for appropriate $k_i \in \Z$.
Since infinite torsion groups are not boundedly generated, it follows that each of the three classes of just-infinite groups contains a finitely generated group that does not have the bounded generation property.
On the other hand, it was proven by Carter and Keller~\cite{CarterKeller83} that $\PSL_n(\Z)$ is boundedly generated for $n \geq 3$, which provides an interesting boundedly generated hereditarily just-infinite group.
The existence of boundedly generated, infinite, simple groups was established by Muranov~\cite{Muranov05}, whose construction seems to be the only one available at present.
It remains to study just-infinite branch groups. The question of existence of boundedly generated just-infinite branch groups was raised by Bartholdi, Grigorchuk, and \v{S}uni\'{k}~\cite[Question 12]{BGS-branch} and remained open to the best of our knowledge.
The purpose of the paper is to show that the answer is negative for arbitrary branch groups (even without the assumption of being just-infinite). 

\begin{theorem}\label{introthm:branchgroups}
There is no boundedly generated branch group.
\end{theorem}

As a consequence, it follows from Wilson's classification of just-infinite groups that every boundedly generated infinite group has a quotient that is virtually a product of finitely many copies of a boundedly generated simple or hereditarily just infinite group.
The proof the Theorem 1 is a rather direct combination of results of Pyber and Segal~\cite{PyberSegal07}, Shusterman~\cite{Shusterman16}, and Ab\'{e}rt~\cite{Abert06}.

 Ab\'{e}rt proved that weakly branch groups are not linear over any field (for 
 branch groups this result is due to Grigorchuk and Delzant). 
More precisely, he defined for every field $k$ the natural number $\mathrm{mat}_k(n)$ to be the minimal $r$ such that every graph on $n$ vertices can be represented in the matrix algebra $\mathrm{M}_{r,r}(k)$ where the graph's edges encode non-commutation.  Ab\'{e}rt showed that $\sqrt{\lfloor n/2 \rfloor} \leq \mathrm{mat}_k(n) \leq 2(n - \lfloor\log_2(n)\rfloor + 1)$ and asked for a linear lower bound \cite[Question 4]{Abert06}. The following result answers this question in the affirmative.

\begin{theorem}\label{introthm:lower-bound-comm}
Let $k$ be a field and let $r \geq 1$.
Suppose that there are $(r\times r)$-matrices $a_1,\dots,a_n, b_1, \dots, b_n \in \mathrm{M}_{r,r}(k)$ such that all pairwise commutators are trivial except for $[a_i,b_i] = a_ib_i-b_ia_i$ for all $i \in \{1,\dots,n\}$.
Then $r \geq n+1$.
\end{theorem}

The lower bound in Theorem \ref{introthm:lower-bound-comm} is sharp. Let $ \lambda \in k^\times$. Consider the matrices $a_i = I + E_{1,i+1}$, $b_i = I - \lambda E_{i+1,i+1} \in \mathrm{M}_{n+1,n+1}(k)$ for $i=1,\ldots,n$, where $I$ is the identity matrix and $E_{i,j}$ denotes the elementary matrix whose $(i,j)$-entry is $1$ and all other entries are $0$.
Then $a_i,b_i$ satisfy the assumptions of Theorem 2. If $\lambda \neq 1$, then $a_i,b_i$ are invertible. If $|k| > 2$, this shows with \cite[Prop.~5]{Abert06} that
\[
	\Bigl\lfloor \frac{n}{2} \Bigr\rfloor + 1 \leq \mathrm{mat}_k(n) \leq n+1.
\]

The non-linearity of weakly branch groups follows since these groups contain infinite products of non-abelian groups.
Let $\mu_k(G)$ denote the minimal dimension of a faithful, finite dimensional representation of a group $G$ over a field $k$ (we write $\mu_k(G) = \infty$ if $G$ is not linear over $k$).
Theorem \ref{introthm:lower-bound-comm} directly implies a lower bound $\mu_k(G_1\times\ldots\times G_n) \geq n+1$ for direct products of non-abelian groups.
Similarly Theorem \ref{introthm:lower-bound-comm}
provides lower bounds for representations of products non-commutative (Lie) algebras.
If the factors $G_i$ are assumed to be non-solvable, the lower bound can be improved further.
\begin{theorem}\label{introthm:non-solvable}
Let $k$ be a field, let $G_1,\dots, G_n$ be groups and let $G = G_1 \times \dots \times G_n$ denote their direct product.
\begin{enumerate}
\item If the groups $G_1,\dots,G_n$ are non-abelian, then $\mu_k(G) \geq n+1$.
\item If the groups $G_1,\dots,G_n$ are non-solvable, then $\mu_k(G) \geq 2n$.
\end{enumerate}
\end{theorem}
Both lower bounds in Theorem \ref{introthm:non-solvable} are sharp.
Let $a_i = I + E_{1,i+1}$, $b_i = I - 2E_{i+1,i+1} \in \GL_{n+1}(\Q)$  be as above.  Setting $G_i = \langle a_i,b_i \rangle$ we can therefore deduce that $\mu_{\Q}(G_1 \times \ldots \times G_n) = n+1$.
Suppose that the groups $G_i$ in Theorem 3 are non-solvable subgroups of $\GL_2(k)$ for some field $k$.
Then each $G_i$ can be embedded in a $2 \times 2$-diagonal block in $\GL_{2n}(k)$, which gives us an embedding of $G = G_1 \times \dots \times G_n$ in $\GL_{2n}(k)$.
Together with Theorem 3 this implies $\mu_k(G) = 2n$. 
In particular, this applies to the case where each $G_i$ is a non-abelian free group and thereby recovers~\cite[Theorem 3]{CampagnoloKammeyer21} in the $\SL_n$-case.

\section{Branch groups are not boundedly generated}

There are several characterizations of branch groups.
The following one, which is a slight reformulation of~\cite[Definition 1.1]{BGS-branch}, does not involve a rooted tree which makes it rather abstract.
However it suits well for our purposes.
A more geometric definition can be found in~\cite[Definition 1.13]{BGS-branch}.

\begin{introdefinition}\label{def:branch-group-abstract}
A group $G$ is called a \emph{branch group} if it admits a decreasing sequence of subgroups $(H_i)_{i \in \N_0}$ with $H_0 = G$ and $\bigcap \limits_{i \in \N_0} H_i = 1$, and a sequence of integers $(k_i)_{i \in \N_0}$ with $k_0 = 1$ such that for each $i$ the following hold:
\begin{enumerate}
\item $H_i$ is a normal subgroup of finite index in $G$.
\item $H_i$ splits as a direct product $H_i = H^{(1)}_i \times \ldots \times H^{(k_i)}_i$, where the factors are pairwise isomorphic.
\item the quotient $m_{i+1} \defeq k_{i+1} / k_i$ is an integer with $m_{i+1} \geq 2$, and the product decomposition of $H_{i+1}$ refines the product decomposition of $H_i$ in the sense that each factor $H^{(j)}_i$ of $H_i$ contains the factors $H^{(\ell)}_{i+1}$ of $H_{i+1}$, where $\ell$ satisfies $(j-1) \cdot m_{i+1}+1 \leq \ell \leq j \cdot m_{i+1}$.
\item $G$ acts transitively by conjugation on the set of factors $H^{(j)}_i$ of $H_i$.
\end{enumerate}
\end{introdefinition}

As indicated in the introduction, not every branch group is just-infinite.
In fact there is no need for finitely generated branch groups to admit a just-infinite quotient that is a branch group.
See~\cite[Theorem 2]{DelzantGrigorchuk08} for an example of finitely generated branch group that maps homomorphically onto $\Z$, which is of course just-infinite and bounded generated.
As a consequence, to prove that branch groups cannot be boundedly generated, it is not sufficient to consider the just-infinite case, in which the claim turns out to be a direct consequence of results of Ab\'{e}rt~\cite{Abert06}, Pyber and Segal~\cite{PyberSegal07}.

\begin{proof}[Proof of Theorem 1]
Suppose there is a branch group $G$ that is boundedly generated.
Then~\cite[Corollary 1.6]{PyberSegal07} tells us that $G$ admits an epimorphism $\pi \colon G \rightarrow Q$, where $Q$ is an infinite linear group.
However, by~\cite[Corollary 7]{Abert06} branch groups are not linear over any field.
Thus $Q$ is a proper quotient of $G$.
As such $Q$ is virtually abelian by~\cite[Proposition 6]{DelzantGrigorchuk08}.
Since $G$, being a boundedly generated group, is finitely generated, it follows that $Q$ has a (non-trivial) free abelian finite index subgroup $Q_0$.
We can therefore consider the finite index subgroup $G_0 \defeq \pi^{-1}(Q_0)$ of $G$, which by construction maps onto $\Z$.
Let us now fix an arbitrary number $n \in \N$.
From the definition of a branch group it follows that $G$ contains a finite index subgroup of the form $H_i = H^{(1)}_i \times \ldots \times H^{(k_i)}_i$, where the factors are pairwise isomorphic and $k_i \geq n$.
Then $H_i \cap G_0$ is a finite index subgroup of $H_i$.
In this case it can be easily seen that there are pairwise isomorphic finite index subgroups $K^{(j)}_i \leq H^{(j)}_i$ such that $K_i \defeq K^{(1)}_i \times \ldots \times K^{(k_i)}_i \leq H_i \cap G_0$.
In particular we see that $K_i$ has finite index in $G_0$, which implies that it maps onto $\Z$.
Thus some, and hence every, factor $K^{(j)}_i$ maps onto $\Z$.
We can therefore deduce that the torsion-free part of the abelianization of $K_i$ has rank at least $k_i \geq n$.
As a consequence, this holds for every finite index subgroup of $K_i$.
In particular this tells us that there is no finite index subgroup of $K_i$ that can be generated with less then $n$ elements.
Since $n \in \N$ was arbitrary this contradicts a result of Shusterman~\cite[Theorem 1.1]{Shusterman16}, which tells us that for every boundedly generated group $H$ there is a constant $C > 0$ such that every finite index subgroup of $H$ contains a finite index subgroup that can be generated by at most $C$ elements.
\end{proof}

\section{Lower bounds for the minimal representation dimension of directs products}

Let us now prove the results concerning the minimal representation dimensions.

\begin{proof}[Proof of Theorem 2]
For the proof we combine ideas from \cite{Abert06} and \cite{LeandroRojas}.
Extending scalars, we may assume that $k$ is an infinite field.
Recall that $I \in \mathrm{M}_{r,r}(k)$ denotes the identity matrix.
We claim that $I, a_1,\dots,a_n,b_1,\dots,b_n$ are linearly independent. This follows along the lines of \cite[Proof of Thm.~3]{Abert06}. Suppose that $cI + \sum_{j} \lambda_j a_j + \sum_{j} \lambda'_j b_j = 0$ for $c,\lambda_1,\dots,\lambda_n,\lambda_1',\dots,\lambda_n' \in k$. Taking commutators with $a_i$ (resp.\ $b_i$) shows $\lambda_i = 0$ (resp.\ $\lambda_i' = 0$); since $I \neq 0$ the last remaining coefficient $c$ vanishes as well.

Let $V = k^r$ and let $C$ denote the linear span of $\{I,a_1,\dots,a_n,b_1,\dots,b_n\}$ in $\mathrm{M}_{r,r}(k)$.
Consider the linear map $\Psi \colon C \to V$ defined by $\Psi(X) = Xv$ for some $v \in V$. We will see that the image of $\Psi$ has dimension at least $n+1$ if $v$ is chosen appropriately.
As the commutators $z_i = [a_i,b_i]$ are non-trivial, the kernel of each $z_i$ is a proper subspace of $V$. However, $V$ cannot be covered by a finite union of proper subspaces (as $k$ is infinite). Thus there is a vector $v \in V$ such that $z_iv \neq 0$ for all $i \in \{1,\dots, n\}$. Let $\alpha \colon V \to k$ be a linear form such that $\alpha(v) \neq 0$ and $\alpha(z_iv) \neq 0$ for all $i \in \{1,2,\dots,n\}$ (such a linear form $\alpha$ exists, as the dual space $V^*$ cannot be covered by finitely many proper subspaces).
Now $\beta \colon C \times C \to k$ defined by $\beta(x,y) = \alpha([x,y](v))$ is an alternating form on $C$. 
It is not difficult to see that $\beta$ is non-degenerate on the subspace $\langle a_1,\dots,a_n,b_1,\dots,b_n\rangle \subseteq C$ (e.g. the matrix representation has full rank).
Let us observe that $kI+\ker(\Psi)$ is an isotropic subspace, since for $x,y \in kI+\ker(\Psi)$ we have $[x,y](v) = xyv-yxv = 0$. 
As $v \neq 0$ we have $I \not\in \ker(\Psi)$ and thus $\dim_k\ker(\Psi) +1 \leq  n+1$. This allows us to conclude that
\[
r \geq \dim_k( \mathrm{im}(\Psi)) = 2n +1 - \dim_k\ker(\Psi) \geq n+1. \qedhere
\]
\end{proof}
\begin{proof}[Proof of Theorem 3]
The first assertion follows immediately from Theorem 2.
Assume now that each $G_i$ is non-solvable.
 If $G$ is not linear, there is nothing to show.
Assume that $(\rho,V)$ is a finite dimensional faithful representation over $k$. By extension of scalars, we may assume that $k$ is algebraically closed.
Let $V^{1},\dots, V^{t}$ denote the composition factors of $V$ considered as $G$-module. Since $k$ is algebraically closed, the composition factor $V^{j}$ is isomorphic to a tensor product
\[
	V^{j} = V_1^{j} \otimes_k V_2^{j} \otimes_k \dots \otimes_k V_n^{j}
\]
where $V^j_i$ is an irreducible $G_i$-representation; see e.g.~\cite[Prop.~2.3.23]{Kowalski}.
The composition factors of $V|_{G_i}$ are the irreducible representations  $V_i^{1},\dots,V_i^t$ each one possibly occurring several times.
Suppose for a contradiction that $V_i^{j}$ is one-dimensional for all $j$.
Then there is a basis of $V$ such that $\rho(G_i)$ is represented by upper triangular matrices. This gives a contradiction, since $G_i$ is not solvable.

For each $j$ let $S_j\subseteq \{1,\dots,n\}$ be the set of $i$ such that $\dim_k V_i^j \geq 2$. 
By the observation above, each $i \leq n$ belongs to at least one of the sets $S_j$. This implies
\[
	 \dim_k V = \sum_{j = 1}^t \prod_{i=1}^n \dim_k V^j_i \geq \sum_{j = 1}^t 2^{|S_j|} \geq \sum_{j=1}^t 2|S_j| \geq 2n. \qedhere\]
\end{proof}

\bibliographystyle{amsplain}
\bibliography{literatur}

\providecommand{\bysame}{\leavevmode\hbox to3em{\hrulefill}\thinspace}
\providecommand{\MR}{\relax\ifhmode\unskip\space\fi MR }
\providecommand{\MRhref}[2]{%
  \href{http://www.ams.org/mathscinet-getitem?mr=#1}{#2}
}
\providecommand{\href}[2]{#2}
\begin{thebibliography}{10}

\bibitem{Abert06}
Mikl\'{o}s Ab\'{e}rt, \emph{Representing graphs by the non-commuting relation},
  Publ. Math. Debrecen \textbf{69} (2006), no.~3, 261--269. \MR{2273978}

\bibitem{BartholdiGrigorchuk99}
Laurent Bartholdi and Rostislav~I. Grigorchuk, \emph{On parabolic subgroups and
  {H}ecke algebras of some fractal groups}, Serdica Math. J. \textbf{28}
  (2002), no.~1, 47--90. \MR{1899368}

\bibitem{BGS-branch}
Laurent Bartholdi, Rostislav~I. Grigorchuk, and Zoran \v{S}uni\'{k},
  \emph{Branch groups}, Handbook of algebra, {V}ol. 3, Handb. Algebr., vol.~3,
  Elsevier/North-Holland, Amsterdam, 2003, pp.~989--1112. \MR{2035113}

\bibitem{BurgerMozes00}
Marc Burger and Shahar Mozes, \emph{Lattices in product of trees}, Inst. Hautes
  \'{E}tudes Sci. Publ. Math. (2000), no.~92, 151--194 (2001). \MR{1839489}

\bibitem{LeandroRojas}
Leandro Cagliero and Nadina Rojas, \emph{Faithful representations of minimal
  dimension of current {H}eisenberg {L}ie algebras}, Internat. J. Math.
  \textbf{20} (2009), no.~11, 1347--1362. \MR{2584190}

\bibitem{CampagnoloKammeyer21}
Caterina Campagnolo and Holger Kammeyer, \emph{Products of free groups in lie
  groups}, Journal of Algebra \textbf{579} (2021), 237--255.

\bibitem{CarterKeller83}
David Carter and Gordon Keller, \emph{Bounded elementary generation of
  $\mathrm{SL}_{n}(\mathcal{O})$}, Amer. J. Math. \textbf{105} (1983), no.~3,
  673--687. \MR{704220}

\bibitem{DelzantGrigorchuk08}
Thomas Delzant and Rostislav Grigorchuk, \emph{Homomorphic images of branch
  groups, and {S}erre's property ({FA})}, Geometry and dynamics of groups and
  spaces, Progr. Math., vol. 265, Birkh\"{a}user, Basel, 2008, pp.~353--375.
  \MR{2402409}

\bibitem{ErshovJaikin10}
Mikhail Ershov and Andrei Jaikin-Zapirain, \emph{Property ({T}) for
  noncommutative universal lattices}, Invent. Math. \textbf{179} (2010), no.~2,
  303--347. \MR{2570119}

\bibitem{Grigorchuk80}
R.~I. Grigor\v{c}uk, \emph{On {B}urnside's problem on periodic groups},
  Funktsional. Anal. i Prilozhen. \textbf{14} (1980), no.~1, 53--54.
  \MR{565099}

\bibitem{Kowalski}
Emmanuel Kowalski, \emph{An introduction to the representation theory of
  groups}, Graduate Studies in Mathematics, vol. 155, American Mathematical
  Society, Providence, RI, 2014. \MR{3236265}

\bibitem{Margulis91}
G.~A. Margulis, \emph{Discrete subgroups of semisimple {L}ie groups},
  Ergebnisse der Mathematik und ihrer Grenzgebiete (3) [Results in Mathematics
  and Related Areas (3)], vol.~17, Springer-Verlag, Berlin, 1991. \MR{1090825}

\bibitem{McCarthy68}
Donald McCarthy, \emph{Infinite groups whose proper quotient groups are finite.
  {I}}, Comm. Pure Appl. Math. \textbf{21} (1968), 545--562. \MR{237637}

\bibitem{Muranov05}
Alexey Muranov, \emph{Diagrams with selection and method for constructing
  boundedly generated and boundedly simple groups}, Comm. Algebra \textbf{33}
  (2005), no.~4, 1217--1258. \MR{2136699}

\bibitem{Olshanskii79}
Alexander~Yu. Ol'shanskii, \emph{Infinite groups with cyclic subgroups}, Dokl.
  Akad. Nauk SSSR \textbf{245} (1979), no.~4, 785--787. \MR{527709}

\bibitem{PyberSegal07}
L\'{a}szl\'{o} Pyber and Dan Segal, \emph{Finitely generated groups with
  polynomial index growth}, J. Reine Angew. Math. \textbf{612} (2007),
  173--211. \MR{2364077}

\bibitem{Shusterman16}
Mark Shusterman, \emph{Ranks of subgroups in boundedly generated groups}, Bull.
  Lond. Math. Soc. \textbf{48} (2016), no.~3, 539--547. \MR{3509913}

\bibitem{Wilson71}
J.~S. Wilson, \emph{Groups with every proper quotient finite}, Proc. Cambridge
  Philos. Soc. \textbf{69} (1971), 373--391. \MR{274575}

\end{thebibliography}

\end{document}